\documentstyle[12pt]{amsart}
\begin{document}
\title{Projectively flat connections and flat connections on homogeneous spaces}
\author{Hajime URAKAWA}
\address{Division of Mathematics\\
Graduate School of Information Sciences\\
Tohoku University\\
Aoba 6-3-09, Sendai, 980-8579, Japan}
\title[Projectively flat connections and flat connections]
{Projectively flat connections and flat connections on homogeneous spaces}
\email{urakawa@@math.is.tohoku.ac.jp}
\thanks{Supported by the Grant-in-Aid for the Scientific Research, (B), 
No. 16340044,  and also 
(A), No. 19204004, (C), No. 21540207, 
Japan Society for the Promotion of Science.}
\dedicatory{Dedicated to the memory of Professor Masaru Takeuchi}
\maketitle
\begin{abstract} 
We show a correspondence between 
the set of all $G$-invariant projectively flat connections 
on a homogeneous space $M=G/K$, 
and 
the one of all $\widetilde{G}$-invariant 
flat connections 
on homogeneous spaces 
$\widetilde{M}=\widetilde{G}/K$, where 
$\widetilde{G}$ is a central extension of $G$ (Theorem 3.3). 
\end{abstract}
\numberwithin{equation}{section}
\theoremstyle{plain}
\newtheorem{df}{Definition}[section]
\newtheorem{th}[df]{Theorem}
\newtheorem{prop}[df]{Proposition}
\newtheorem{lem}[df]{Lemma}
\newtheorem{cor}[df]{Corollary}
\newtheorem{rem}[df]{Remark}
\section{Introduction and statement of results.} 
Flat connections and projectively flat ones 
have been extensively studied 
by many authors (for examples, 
\cite{A}, \cite{D}, \cite{E}, \cite{MO}, \cite{N}, \cite{NS}, \cite{S2}, \cite{U}).  
Even though they are of course very different objects,  
but it seems that there would exist deep unknown relations between 
them each other. 
In this paper, we want to show 
some relation between 
the set of all $G$-invariant projectively flat connections 
on a homogeneous space $M=G/K$, 
and 
the one of all $\widetilde{G}$-invariant 
flat connections 
on homogeneous spaces 
$\widetilde{M}=\widetilde{G}/K$, where 
$\widetilde{G}$ is a central extension of $G$. 
\par
Indeed, let $\widetilde{G}\supset G\supset K$ be three Lie groups 
with Lie algebras 
$\widetilde{\frak g}\supset {\frak g}\supset {\frak k}$, 
where $\widetilde{\frak g}={\frak g}\oplus{\Bbb R}E$, 
and $[E,\widetilde{\frak g}]=\{0\}$. 
Let us consider two homogeneous (not necessarily reductive) 
spaces 
$\widetilde{M}=\widetilde{G}/K$ and 
$M=G/K$, respectively. 
We consider the two sets of all $\widetilde{G}$-invariant 
{\em flat} affine connections on $\widetilde{M}$, and 
of all $G$-invariant {\em projectively flat}
affine connections on $M$, 
which correspond to 
the sets ${\mathcal F}_0(\widetilde{\frak g},{\frak k})$ and 
${\mathcal PF}_0(\widetilde{\frak g},{\frak k})$
of irreducible affine representations of $\widetilde{\frak g}$, 
respectively (for more precise, see \S 3, Definitions 3.1, 3.2). 
Then, 
\begin{th} $($cf. Theorem 3.3$)$
It holds that 
$${\mathcal F}_0(\widetilde{\frak g},{\frak k})=
{\mathcal PF}_0(\widetilde{\frak g},{\frak k})\cup 
{\mathcal F}_0^{I\!I}(\widetilde{\frak g},{\frak k}),$$
where
${\mathcal F}_0^{I\!I}(\widetilde{\frak g},{\frak k})$
is 
the set of all real irreducible affine representations 
$(\widetilde{f},\widetilde{q}, \widetilde{V})$ of 
$\widetilde{\frak g}$ satisfying that 
\par\qquad
$(1)$  
$\dim \widetilde{V}=\dim G/K+1$,
\par\qquad 
$(2)$ $\widetilde{V}$ admits  an 
$\widetilde{f}(\widetilde{\frak g})$-invariant complex structure $J$, 
and 
\par\qquad $(3)$ 
there exists a non-zero element $v_0\in \widetilde{V}$ satisfying that 
$$
\widetilde{f}({\frak k})v_0=\{0\}, \,\,\text{and}\,\,\,\,
\widetilde{V}=\widetilde{f}({\frak g})v_0
\oplus {\Bbb R}\widetilde{f}(E)v_0.
$$
\par
In particular, in the case that 
$\dim M=\dim G/K$ is even, then, 
$${\mathcal F}_0(\widetilde{\frak g},{\frak k})=
{\mathcal PF}_0(\widetilde{\frak g},{\frak k})
\quad \text{and}\quad 
{\mathcal F}^{I\!I}_0(\widetilde{\frak g},{\frak k})=\emptyset.$$
Thus, $G/K$ admits a $G$-invariant flat connection if and only if $\widetilde{G}/K$ admits a $\widetilde{G}$-invariant projectively flat connection in the case that $G/K$ is of even dimension. 
\end{th}
\vskip0.6cm
\par
Let us recall a classification of real semi-simple Lie groups admitting a projectively flat connection $($\cite{A}, \cite{U}, \cite{E}$)$.  
\begin{th}
Let $G$ be a real semi-simple Lie group. Then, $G$ admits a left invariant projectively flat connection if and only if the Lie algebra 
$\frak g$ is one of the following:
\vskip0.6cm
\par
$(1) \quad {\frak sl}(n+1,{\Bbb R})$, \quad $n\geq 1$,
\par
$(2) \quad {\frak su}^{\ast}(2n)$, \quad $n\geq 1$, 
\par
\noindent
where ${\frak su}^{\ast}(2n)$ is 
the Lie algebra given by 
$$
{\frak su}^{\ast}(2n)=\left\{
\begin{pmatrix}
Z_1 & Z_2\\
-\overline{Z}_2 & \overline{Z}_1
\end{pmatrix}
;\,Z_1,\,Z_2\in M(n,{\Bbb C}), \text{\rm Tr}Z_1+\text{\rm Tr}
\overline{Z}_1=0
\right\}. 
$$
\end{th}
\vskip0.6cm\par
Since Theorems 1.1 and 1.2 (and also Remark 1.6) 
except the case of the real representation of 
$G=SL(n+1,{\mathbb R})$ on 
${\frak gl}(n+1,{\mathbb R})$ 
(cf. Theorem 4.3 in Chapter 3 in \cite{NS}, see also \cite{A}), we have 
\begin{cor} 
Let $G$ be a real semi-simple Lie group with Lie algebra $\frak g$ 
of {\em even dimension}, and let 
$\widetilde{\frak g}={\frak g}\oplus {\Bbb R}E$, with $[E,\widetilde{\frak g}]=\{0\}$. 
Let $\widetilde{G}$ be a simply connected Lie group with Lie algebra 
$\widetilde{\frak g}$. 
Then, $\widetilde{G}$ admits a left invariant {\em flat} affine connection 
if and only if $G$ admits a left invariant {\em projectively flat} affine connection. 
In this case,  
$\frak g$ is ${\frak sl}(n+1,{\Bbb R})$,  
where $n\geq 1$ is even. 
\end{cor}
\vskip0.6cm
\par
Let us recall a classification of irreducible Riemannian symmetric spaces admitting invariant projectively flat connections (\cite{A}, \cite{U}). 
\begin{th} Let $M=G/K$ be an irreducible simply connected Riemannian symmetric space. Then, $M$ admits a $G$-invariant {\em projectively flat} affine connection if and only if $M=G/K$ is one of the following: 
\vskip0.6cm
\par
$(1)$ $S^n=SO(n+1)/SO(n)$ \quad $n\geq 2$,
\par
$(2)$ $SL(n+1,{\Bbb R})/SO(n+1)$ \quad $n\geq 2$,
\par
$(3)$ $SU^{\ast}(2n)/Sp(n)$ \quad $n\geq 3$,
\par
$(4)$ $SO_0(n,1)/SO(n)$ \quad $n\geq 2$,
\par
$(5)$ $SL(n+1,{\Bbb C})/SU(n+1)$ \quad $n\geq 1$
\par
$(6)$ $E_6/F_4$ $($of non-compact type EIV$)$.
\end{th}
\vskip0.6cm\par
Since Theorems 1.1 and 1.4 (and also Remark 1.6), we have
\begin{cor} Let $M=G/K$ be an irreducible simply connected Riemannian symmetric space of 
{\em even dimension}, $\widetilde{\frak g}={\frak g}\oplus {\Bbb R}E$ 
with $[E,\widetilde{\frak g}]=\{0\}$, and $\widetilde{G}$, a Lie group with Lie algebra $\widetilde{\frak g}$ 
and $\widetilde{M}$ is simply connected. Then, 
$M=G/K$ admits a $G$-invariant {\em projectively flat} affine connection if and only if $\widetilde{M}=\widetilde{G}/K$ admits a $\widetilde{G}$-invariant 
{\em flat} affine connection. In this case, $G/K$ is one of $(1)\sim (6)$ in Theorem 1.4 of even dimension. 
\end{cor}
\vskip0.6cm\par         
\begin{rem} 
In Corollaries 1.3, and 1.5, Agaoka $($\cite{A}$)$ showed that, if $G$ $($resp. $M=G/K$$)$ admits a left invariant $($resp. invariant$)$ flat connection, then $\widetilde{G}$ 
(resp. $\widetilde{M}=\widetilde{G}/K$) admits a left invariant 
$($resp. $\widetilde{G}$-invariant$)$ flat connection. In this, paper, we show the reverse of his results 
in Corollaries 1.3 and 1.5 in the case of even dimension. 
\end{rem}
\vskip0.6cm\par
{\bf Acknowledgement}: The author would like to express his gratitude to Professor T. Yamada 
who gave useful discussions during his stay 
at Tohoku University as his postdoctoral student from April, 2006, 
to September, 2007, and discussions at January, 2009, at Shimane University, and Professor A. Kasue for the financial support, 
and the referee 
who pointed out  several errors in the manuscript and gave many suggestions 
to improve it into this form.  
\vskip0.6cm\par
\section{Preliminaries.}
In this section, we prepare materials and several facts on
invariant connection on homogeneous spaces (cf. \cite{KN}) and also 
invariant flat connections and projectively flat invariant connections on 
homogeneous spaces (cf. \cite {S1}). 
\par
Let us consider a $C^{\infty}$ affine connection $D$ on a $C^{\infty}$ manifold $M$. In this paper, we assume that everything is $C^{\infty}$. 
\begin{df}
A connection 
$D$ is to be {\em flat} if it has a vanishing  curvature tensor 
$R^D$ and a vanishing torsion tensor $T^D$, 
where 
$$
R^D(X,Y)=D_XD_Y-D_YD_X-D_{[X,Y]},
$$
and 
$$
T^D(X,Y)=D_XY-D_YX-[X,Y]
$$
for all vector fields $X$, $Y$ on $M$. 
\end{df}
\vskip0.6cm\par
Let us recall also that the theory of invariant connections on homogeneous spaces 
(cf. \cite{KN}, Vol. II, p.188). In case of affine connections, we have 
\begin{th}
Let $M$ be a homogeneous space $M=G/K$. Then, there exists a one-to-one correspondence 
between the set of $G$-invariant affine connections on $M=G/K$ and the set of linear mappings 
$\Lambda:\,{\frak g}\rightarrow {\frak gl}(T_o(G/K))$ $($$o=\{K\}$, the origin of $G/K$$)$ such that 
\begin{align}
\Lambda(X)&=\lambda(X) \quad (X\in {\frak k}),\\
\Lambda ({\rm Ad}(k)(X))&={\rm Ad}(\lambda(k))(\Lambda(X))\quad (k\in K, X\in {\frak g}), 
\end{align}
where $\lambda$ is the isotropy representation of $K$, i.e., for $k\in K$, 
$\lambda(k)=k_{\ast}:\,T_o(G/K)\rightarrow T_o(G/K)$ denotes the differential of $k$ at 
$o=\{K\}$. 
\par
To each $G$-invariant connection $D$ on $M=G/K$, there corresponds the linear mapping 
$\Lambda$ defined by 
\begin{equation}
\Lambda(X)=-(A_{\widetilde{X}})_o\quad (X\in {\frak g}).
\end{equation}
Here,  each $X\in {\frak g}$ induces a tangent vector $X_o\in T_oM$, and 
also a vector field $\widetilde{X}$ on $M=G/K$, naturally. 
For all vector field $\widetilde{X}$ on $M$, $A_{\widetilde{X}}$ is the tensor field 
of type (1,1) on $M$ defined by
\begin{equation}
A_{\widetilde{X}}=L_{\widetilde{X}}-D_{\widetilde{X}},
\end{equation}
where
$L_{\widetilde{X}}$ is Lie derivative by $\widetilde{X}$. 
\end{th}
\vskip0.6cm\par
Assume that a homogeneous space $M=G/K$ is {\it reductive}, i.e., 
the Lie algebra ${\frak g}$ is decomposed into ${\frak g}={\frak k}\oplus {\frak m}$, 
where ${\frak k}$ is the Lie algebra of $K$ and ${\frak m}$ is an 
Ad$(K)$-invariant subspace of ${\frak g}$. 
Then, we have (\cite{KN}, Vol II, p. 191) 
\begin{th}
Assume that $M=G/K$ is a reductive homogeneous space with decomposition 
${\frak g}={\frak k}\oplus {\frak m}$. Then, there is a one-to-one correspondence between 
the set of $G$-invariant affine connections on $M$ and the set of linear mappings 
$\Lambda_{\frak m}:\,{\frak m}\rightarrow {\frak gl}(T_o(G/K))$ such that 
\begin{equation}
\Lambda_{\frak m}({\rm Ad}(k)(X))={\rm Ad}(\lambda(k))(\Lambda_{\frak m}(X)) 
\quad (X\in {\frak m},\,k\in K),
\end{equation}
where $\lambda(k)$ denotes the isotropy representation of $K$ on $G/K$. 
The correspondence is given by 
\begin{equation}
\Lambda(X)=
\left\{
\begin{aligned}
&\lambda(X)\qquad (X\in {\frak k}),\\
&\Lambda_{\frak m}(X) \quad (X\in {\frak m}). 
\end{aligned}
\right.
\end{equation}
To each $G$-invariant connection $D$,  the correspondence (2.6) is 
given by 
\begin{equation}
\Lambda_{\frak m}(X)=-(A_{\widetilde{X}})_o\quad (X\in {\frak m}).
\end{equation}
\end{th}
The torsion tensor $T^D$ and the curvature tensor $R^D$ of a $G$-invariant connection 
$D$ can be expressed in terms of $\Lambda$ as follows (\cite{KN}, Vol. II, p. 189).
\begin{th}
In the above theorems, 
the torsion tensor $T^D$ and the curvature tensor $R^D$ of 
a $G$-invariant connection $D$ can be expressed as follows: 
\begin{align}
T^D(X,Y)_o&=
\Lambda(X)(Y_o)-\Lambda(Y)(X_o)-[X,Y]_o\quad (X,\,Y\in {\frak g}),\\
R^D(X,Y)_o&=
[\Lambda(X),\Lambda(Y)]-\Lambda([X,Y])
\quad (X,\,Y\in {\frak g}). 
\end{align}
\end{th}
\vskip0.6cm\par
For $G$-invariant flat connections on $M=G/K$, 
due to the above theorems, we have
(\cite{S1}) 
\begin{th}
Assume that a homogeneous space $M=G/K$ admits a $G$-invariant flat connection on $M=G/K$. 
Then, there exists a affine representation
$(f,q,V)$ of $\frak g$ on $V$ such that 
\begin{equation}
\left\{
\begin{aligned}
&\dim V=\dim M,\\
&q :\,{\frak g}\rightarrow V\,\,\text{\rm is surjective, and \,\,Ker}(q)={\frak k}.
\end{aligned}
\right.
\end{equation}
Conversely, if $G$ is simply connected, and the Lie algebra $\frak g$ 
admits an affine representation $(f,q,V)$ satisfying (2.10), then $M=G/K$ admits a $G$-invariant flat affine connection. 
Here, two linear mappings $q:\, {\frak g}\rightarrow V$, and 
$f:\,{\frak g}\rightarrow {\frak gl}(V)$ are given by 
$q(X)=X_o$ $(X\in {\frak g})$, and 
$f(X)=\Lambda(X)$ $(X\in {\frak g})$, where $V=T_oM$.
Then, that $(f,q,V)$ is an affine representation 
${\frak g}$ on $V$ means that 
\begin{equation}
\left\{
\begin{aligned}
&[f(X),f(Y)]=f([X,Y])\quad (X,Y\in {\frak g}),\\
&q([X,Y])=f(X)q(Y)-f(Y)q(X)\quad (X,Y\in {\frak g}).
\end{aligned}
\right.
\end{equation}
\end{th}
\vskip0.6cm\par
Now, let us recall the notion of projectively flat connections. 
\begin{df}
$D$ is to be {\em projectively flat} if the Ricci tensor 
$$
\text{\rm Ric}^D(Y,Z):=
\text{\rm Tr}(\{T_pM\ni X\mapsto R^D(X,Y)Z\in T_pM\})
$$
$(Y,Z\in T_pM$ $p\in M)$, 
is {\em symmetric}, i.e., $\text{\rm Ric}^D(Y,Z)=\text{\rm Ric}^D(Z,Y)$, and 
for every $p\in M$, there exists a neighborhood $U$ of $p$ such that $D$ is {\em equivalent} to a flat connection $\overline{D}$ on $U$, i.e., there exists a closed 1-form $\rho$ on $U$ such that 
$$
D_XY=\overline{D}_XY+\rho(X)Y+\rho(Y)X,
\quad \text{for all} \,\, X,Y\in {\frak X}(U).
$$
\end{df}
\vskip0.6cm\par
A classical theorem says that 
\begin{th}
Assume that $D$ is an affine connection of which $\text{\rm Ric}^D$ is symmetric and $T^D=0$. Then, $D$ is projectively flat if and only if the Weyl curvature tensor vanishes, i.e., 
$$
R^D(X,Y)Z=\frac{1}{n-1}\{\text{\rm Ric}^D(Y,Z)X-\text{\rm Ric}^D(X,Z)Y\}
$$
and the {\em Codazzi equation} holds, i.e., 
$$
(D_X\text{\rm Ric}^D)(Y,Z)=(D_Y\text{\rm Ric}^D)(X,Z)
$$
for all $X,Y,Z\in {\frak X}(M)$. 
\end{th}
\vskip0.6cm
\par
Let us consider 
a {\em centro-affine} immersion of an $n$-dimensional manifold $M$ into ${\Bbb R}^{n+1}$, $\varphi:\,M \hookrightarrow {\Bbb R}^{n+1}$, that is, 
$$
T_{\varphi(p)}{\Bbb R}^{n+1}=\varphi_{\ast}(T_pM)\oplus {\Bbb R}\,\overrightarrow{o\,\varphi(p)},\quad (p\in M). 
$$
Then, the induced connection $D$ on $M$ from the standard flat connection $D^0$ on ${\Bbb R}^{n+1}$ via $\varphi$, i.e., 
$$
(D^0_XY)_p=\varphi_{\ast}(D_XY)_p+h(X,Y)(-\overrightarrow{o\,\varphi(p)}),\quad (X,Y\in {\frak X}(M),\,p\in M),
$$
is projectively flat. Furthermore, it holds (\cite{S2}) that 
\begin{th}
Let $M=G/K$ be a simply connected homogeneous space. Then, the following two conditions are equivalent: 
\par
$(1)$ \quad $M=G/K$ admits a $G$-invariant {\em projectively flat} connection.
\par
$(2)$\quad $M=G/K$ admits a $G$-equivariant centro-affine immersion. 
\end{th}
\vskip0.6cm
\par
Furthermore, Shima showed (\cite{S2}, see also \cite{S1}, p. 228) that 
\begin{th}
Let $M=G/K$ be an arbitrary homogeneous space. 
\par
$(1)$  Assume that $M=G/K$ admits a $G$-invariant {\em projectively flat} connection. Let $\frak g$ be the Lie algebra of $G$, and $\frak k$ the Lie subalgebra corresponding to $K$, respectively. Let $\widetilde{\frak g}$ be the central extension of $\frak g$, i.e., 
$\widetilde{\frak g}={\frak g}\oplus {\Bbb R}E$, where 
$[E,\widetilde{\frak g}]=\{0\}$. Then, $\widetilde{\frak g}$ admits an {\em affine representation} 
$(\widetilde{f},\widetilde{q},\widetilde{V})$ on a vector space $\widetilde{V}$ of dimension $\dim M+1$ satisfying the following two conditions: 
\par
\qquad $(i)$ \quad $\widetilde{q}:\,\widetilde{\frak g}\rightarrow\widetilde{V}$ is surjective and $\text{\rm Ker}(\widetilde{q})$ is $\frak k$, 
\par
\qquad $(ii)$ \quad $\widetilde{f}(E)$ is the identity map of $\widetilde{V}$ and $\widetilde{q}(E)\not=0$. 
\par
$(2)$ Conversely, if $\widetilde{\frak g}$ admits an affine representation $(\widetilde{f},\widetilde{q},\widetilde{V})$ on $\widetilde{V}$ of dimension $\dim M+1$ 
satisfying $(i)$ and $(ii)$, then, $M=G/K$ admits a $G$-invariant {\em projectively flat} affine connection if $G$ is simply connected. 
\end{th}
\vskip0.6cm
\par
\section{Invariant projectively flat connections and flat connections.}
Let us begin the following example due to Shima 
(\cite{S1}, Ch. 9, Exercise 9.1.1): 
\vskip0.6cm\par
{\bf Example 3.1.} \quad Let $\widetilde{G}=GL(n,{\Bbb R})$, 
$$
K:=\left\{
\begin{pmatrix}
I_r& x\\ O&y
\end{pmatrix}
;\,x\in M(r,n-r), y\in M(n-r,n-r)
\right\}. 
$$
Then, the quotient space $\widetilde{M}=\widetilde{G}/K$ 
admits a $\widetilde{G}$-invariant flat connection. 
Here, $I_r$ is the identity matrix of degree $r$, and 
$O$ is the $(n-r)\times r$ zero matrix where $r=1,\dots,n-1$. 
Then, $\widetilde{\frak g}={\frak gl}(n,{\Bbb R})$, 
$$
{\frak k}:=\left\{
\begin{pmatrix}
O& X\\ O&Y
\end{pmatrix}
;\,X\in M(r,n-r), Y\in M(n-r,n-r)
\right\}, 
$$
and $T_0\widetilde{M}$ is isomorphic to 
$V$
where 
$$
V=
\left\{
\begin{pmatrix}
X'& O\\ Y'&O
\end{pmatrix}
;X'\in M(r,r),Y'\in M(n-r,r)
\right\}.
$$
Then, $\widetilde{\frak g}={\frak k}\oplus V$, and 
for every $X\in \widetilde{\frak g}$, 
$f(X)\in {\frak gl}(V)$ and $q(X)\in V$ are defined by 
$$
f(X)v=Xv\in V, \quad 
q(X)=X_V\in V, \quad (X\in \widetilde{\frak g},\,v\in V),
$$
where $Xv$ is the matrix multiplication of $X$ and $v$, 
and $X_V\in V$ is the $V$-component of $X\in \widetilde{\frak g}$ corresponding to the decomposition of 
$\widetilde{\frak g}={\frak k}\oplus V$. Then, 
$(f,q,V)$ is an affine representation of 
$\widetilde{\frak g}$ on 
$V$ which induces the $\widetilde{G}$-invariant flat connection on 
$\widetilde{M}$.  
But, 
$\widetilde{M}=\widetilde{G}/K$ is not a reductive homogeneous space since 
there is no ad$({\frak k})$-invariant decomposition 
of $\widetilde{\frak g}$ (see also \cite{KN}, Vol. II, p. 199, Example 2.1). 
The action of ad$(\widetilde{\frak g})$ on $V$ is
irreducible if and only if $r=1$. 
\vskip0.6cm\par
{\bf Example 3.2.} \,\, 
$GL(n,{\Bbb R})$, $GL(n,{\Bbb C})$,  
the upper triangular nilpotent Lie group $N_1$: 
$$
N_1=\left\{
\begin{pmatrix}
1&   &&\ast\\
  &1&&\\
  &&\ddots &&\\
  O &&&1
\end{pmatrix}; \,\ast\,\,\text{is arbitrary} 
\right\}, 
$$
the $(2n+1)$-dimensional Heisenberg nilpotent Lie group $N_2$, 
and the solvable Lie group $S$: 
$$
S=\left\{
\begin{pmatrix}
e^{\theta_1}&   &&\ast\\
  &e^{\theta_2}&&\\
  &&\ddots &&\\
  O &&&e^{\theta_n}
\end{pmatrix};\,\theta_1,\theta_2,\dots,\theta_n\in {\Bbb R}, \, 
\ast\,\,\text{is arbitrary} 
\right\}
$$
admit left invariant flat connections. 
\par
But, $SU(2)$ admits no left invariant flat connection. 
\vskip0.6cm\par
In this section, we do not assume that 
every  homogeneous space is reductive.  
Our setting is as follows. 
Let $\widetilde{G}\supset G\supset K$ be three Lie groups with Lie algebras 
$\widetilde{\frak g}\supset {\frak g}\supset {\frak k}$, respectively.  
Let us consider two quotient spaces 
$\widetilde{M}=\widetilde{G}/K$, and $M=G/K$, respectively. 
Assume that 
the Lie algebra $\widetilde{\frak g}$ 
is a central extension of ${\frak g}$ given by 
$\widetilde{\frak g}={\frak g}\oplus {\Bbb R}E$, 
where $[E,\widetilde{\frak g}]=\{0\}$. Then, 
$\dim\widetilde{M}=\dim M+1$. 
\begin{df}
Let us denote by 
${\mathcal F}_{\widetilde{G}}(\widetilde{M})$ 
the set of all $\widetilde{G}$-invariant flat affine connections on 
$\widetilde{M}$, and 
${\mathcal PF}_G(M)$, 
the set of all $G$-invariant projectively flat affine connections 
on $M$. 
Let us define 
${\mathcal F}(\widetilde{\frak g},{\frak k})$, 
the set of all affine representations 
$(\widetilde{f},\widetilde{q},\widetilde{V})$ of 
$\widetilde{\frak g}$ on $\widetilde{V}$ satisfying that 
\par\qquad
$(1)$ $\dim \widetilde{V}=\dim \widetilde{M}$, 
\par\qquad
$(2)$ $\widetilde{q}:\,\widetilde{\frak g}\rightarrow \widetilde{V}$ is a 
surjection, and  {\rm Ker}$(\widetilde{q})={\frak k}$, 
\par\noindent
and also 
${\mathcal PF}(\widetilde{\frak g},{\frak k})$,  
the set of all affine representations 
$(\widetilde{f},\widetilde{q},\widetilde{V})$ of 
$\widetilde{\frak g}$ on $\widetilde{V}$ satisfying that 
\par\qquad
$(1)$ $\dim \widetilde{V}=\dim M+1$, 
\par\qquad
$(2)$ $\widetilde{q}:\,\widetilde{\frak g}\rightarrow \widetilde{V}$ is a 
surjection, 
 {\rm Ker}$(\widetilde{q})={\frak k}$, 
\par\qquad
$(3)$ $\widetilde{f}(E)$ is a non-zero constant multiple of 
the identity 
\par \qquad\quad transformation of $\widetilde{V}$, 
and $\widetilde{q}(E)\not=0$, 
\par\noindent
respectively.  
\end{df} 
Remark here that in the original definition of 
${\mathcal PF}(\widetilde{\frak g},{\frak k})$ 
in \cite{S1} corresponding to the above definition $(3)$
was that:
\par\qquad
$(3')$ {\em $\widetilde{f}(E)$ is the identity transformation of 
$\widetilde{V}$, and $\widetilde{q}(E)\not=0$. }
\par\noindent 
But, if we replace $E$ into a constant multiple of $E$, then we have 
$(3')$ from $(3')$, so that our 
${\mathcal PF}(\widetilde{\frak g},{\frak k})$ corresponds bijectively to 
${\mathcal PF}_G(M)$.
\vskip0.6cm\par 
If $\widetilde{G}$ is simply connected, 
there exists a one-to-one correspondence between 
${\mathcal F}_{\widetilde{G}}(\widetilde{M})$ and 
${\mathcal F}(\widetilde{\frak g},{\frak k})$, and also 
if $G$ is simply connected, 
there exists a one-to-one correspondence between 
${\mathcal PF}_{G}(M)$ and 
${\mathcal PF}(\widetilde{\frak g},{\frak k})$ 
(\cite{S1}). 
By definition, 
${\mathcal PF}(\widetilde{\frak g},{\frak k})
\subset {\mathcal F}(\widetilde{\frak g},{\frak k})$, 
so that 
${\mathcal PF}_G(M)\subset {\mathcal F}_{\widetilde{G}}(\widetilde{M})$. 
\par
To analyze them further, let us define
\begin{df}\quad  
$$
{\mathcal F}_0(\widetilde{g},{\frak k})=
\{(\widetilde{f},\widetilde{q},\widetilde{V})\in {\mathcal F}(\widetilde{\frak g},{\frak k});\,
(\widetilde{f},\widetilde{q},\widetilde{V})\,\,\text{is irreducible under} 
\,\,\widetilde{\frak g}\}, 
$$
and also 
$$
{\mathcal PF}_0(\widetilde{g},{\frak k})=
\{(\widetilde{f},\widetilde{q},\widetilde{V})\in {\mathcal PF}(\widetilde{\frak g},{\frak k});\,(\widetilde{f},\widetilde{q},\widetilde{V})\,\,\text{is irreducible under} 
\,\,\widetilde{\frak g}\}, 
$$
respectively. Then, 
${\mathcal PF}_0(\widetilde{g},{\frak k})\subset
{\mathcal F}_0(\widetilde{g},{\frak k})$ by definition. 
\end{df}
\vskip0.6cm\par
Then, we obtain 
\begin{th} It holds that 
$${\mathcal F}_0(\widetilde{\frak g},{\frak k})=
{\mathcal PF}_0(\widetilde{\frak g},{\frak k})\cup 
{\mathcal F}_0^{I\!I}(\widetilde{\frak g},{\frak k}),$$
where the set 
${\mathcal F}_0^{I\!I}(\widetilde{\frak g},{\frak k})$
is a subset of ${\mathcal F}_0(\widetilde{\frak g},{\frak k})$, 
and coincides with  
the one of all real irreducible affine representations 
$(\widetilde{f},\widetilde{q},\widetilde{V})$ of 
$\widetilde{\frak g}$ satisfying that 
\par\qquad
$(1)$  
$\dim \widetilde{V}=\dim G/K+1$,
\par\qquad 
$(2)$ $\widetilde{V}$ admits  an 
$\widetilde{f}(\widetilde{\frak g})$-invariant complex structure $J$, 
and 
\par\qquad $(3)$ 
there exists a non-zero element $v_0\in \widetilde{V}$ satisfying that 
$$
\widetilde{f}({\frak k})v_0=\{0\}, \,\,\text{and}\,\,\,\,
\widetilde{V}=\widetilde{f}({\frak g})v_0
\oplus {\Bbb R}\widetilde{f}(E)v_0.
$$
\par
In particular, in the case that 
$\dim M=\dim G/K$ is even, then, 
$${\mathcal F}_0(\widetilde{\frak g},{\frak k})=
{\mathcal PF}_0(\widetilde{\frak g},{\frak k})\quad 
\text{and}\quad 
{\mathcal F}^{I\!I}_0(\widetilde{\frak g},{\frak k})=\emptyset.$$
\end{th}
\begin{pf}
\quad 
Let $(\widetilde{f},\widetilde{q},\widetilde{V})
\in {\mathcal F}_0(\widetilde{\frak g},{\frak k})$. 
Then, the affine representation 
$(\widetilde{f},\widetilde{q},\widetilde{V})$ satisfies that 
\par\qquad 
$(1)$ \quad $(\widetilde{f},\widetilde{q},\widetilde{V})$ is an irreducible representation of  $\widetilde{\frak g}$, 
\par\qquad 
$(2)$ \quad $\dim \widetilde{V}=\dim \widetilde{M}=\dim M+1$, 
\par\qquad 
$(3)$ \quad $\widetilde{q}:\,\widetilde{\frak g}\rightarrow \widetilde{V}$ is a surjection, and 
Ker$(\widetilde{q})={\frak k}$. 
\par
By $(1)$, $\widetilde{f}(E)$ is a semi-simple linear transformation of $\widetilde{V}$ because 
$(\widetilde{f},\widetilde{q},\widetilde{V})$ is a completely reducible representation of ${\frak g}$ (see, for example, \cite{T}, p. 28, Theorem 2.9). Then, $\widetilde{V}$ is decomposed into 
$$
\widetilde{V}=V_1\oplus\cdots\oplus V_s\oplus V_{s+1}\oplus\cdots\oplus V_r,
$$ 
where 
\par\quad 
$(i)$ each $V_i$ is $\widetilde{f}(E)$-invariant and the only $\widetilde{f}(E)$-invariant subspaces of $V_i$ are $V_i$ itself or $\{0\}$, and 
\par\quad 
$(ii)$  for the complex extension $\widetilde{f}(E)^{\Bbb C}$ of 
$\widetilde{f}(E)$ to $V_i{}^{\Bbb C}$, 
there exist linearly independent vectors $\{v^i_j\}_{j=1}^{d_i}$ in 
the complexification $V_i{}^{\Bbb C}$ of $V_i$ and 
some complex number $\lambda^i_j$ such that 
\begin{equation}
\widetilde{f}(E)^{\Bbb C}v^i_j=\lambda^i_jv^i_j\quad (j=1,\dots,d_i).
\end{equation}
Notice here that $\dim V_i{}^{\mathbb C}=1$, and $d_i=1$ in this case due to $(i)$. 
\par\quad
$(iii)$ For each $i=1,\dots,s$, $\lambda^i_j=a_i$ is a real number, and $V_i$ itself is the eigenspace of $\widetilde{f}(E)$, and $\{v^i_j\}_{j=1}^{d_i}$ is a basis of $V_i{}^{\mathbb C}$. 
\par\quad 
$(iv)$ For each $i=s+1,\dots,r$, $\lambda^i_j=a^i_j+\sqrt{-1}b^i_j$ 
where $a^i_j$ and $b^i_j$ are real numbers 
with $b^i_j\not=0$, and 
$v^i_j=u^i_j+\sqrt{-1}w^i_j$ $(u^i_j,w^i_j\in V_i;\,j=1,\dots,d_i)$ such that 
\begin{equation}
\left\{
\begin{aligned}
\widetilde{f}(E)u^i_j&=a^i_ju^i_j-b^i_jw^i_j\\
\widetilde{f}(E)w^i_j&=b^i_ju^i_j+a^i_jw^i_j. 
\end{aligned}
\right.
\end{equation}
\par
Then, by $(i)$, it holds that 
$\dim V_i{}^{\mathbb C}=2$, $d_i=1$, and for $j=1$, 
\begin{equation}
\left\{
\begin{aligned}
&\widetilde{f}(E)^{\mathbb C}\overline{v^i_1}
=\overline{\lambda^i_1}\,\overline{v^i_1},\\
&{V_i}^{\mathbb C}=
{\mathbb C}\,v^i_1\oplus{\mathbb C}\overline{v^i_1},
\end{aligned}
\right.
\end{equation}
where $\{v^i_1,\overline{v^i_1}\}$ 
is a basis of $V_i{}^{\mathbb C}$. 
\vskip0.3cm\par
Then we have two cases: $V_1\oplus\cdots\oplus V_s$ is $\{0\}$ or not.  
\vskip0.3cm\par
\underline{Case $(a)$}: $V_1\oplus\cdots\oplus V_s\not=\{0\}$.
\par
In this case, 
$V_1\oplus\cdots\oplus V_s$ is a non-zero 
$\widetilde{f}(\widetilde{\frak g})$-invariant subspace. 
Indeed, for each 
$v\in V_i$ and $\widetilde{X}\in \widetilde{\frak g}$, 
since $[E,\widetilde{\frak g}]=\{0\}$ 
and $\widetilde{f}$ is a Lie algebra homomorphism, we have
$$
\widetilde{f}(E)(\widetilde{f}(\widetilde{X})v)
=\widetilde{f}(\widetilde{X})(\widetilde{f}(E)v)
=a_i\widetilde{f}(\widetilde{X})v,
$$
which implies $\widetilde{f}(\widetilde{X})v$ belongs to 
$V_1\oplus\cdots\oplus V_s$. 
Thus, 
$\widetilde{V}=V_1\oplus\cdots\oplus V_s$ and 
$\widetilde{f}(E)$ has a unique real eigenvalue, say, $a\in {\Bbb R}$,  
since $(\widetilde{f},\widetilde{V})$ is an irreducible representation of 
$\widetilde{\frak g}$. 
\par
Furthermore, $a\not=0$. 
Because if we assume that $a=0$, then, 
$\widetilde{f}(E)=0$. 
Then, for each $\widetilde{X}\in \widetilde{\frak g}$, 
we have 
\begin{equation}
\widetilde{f}(\widetilde{X})(\widetilde{q}(E))=
\widetilde{f}(E)(\widetilde{q}(\widetilde{X}))=0,
\end{equation}
which implies that $\widetilde{q}(E)=0$. 
Because 
if we assume 
$\widetilde{q}(E)\not=0$, 
then, $\{0\}\not={\Bbb R}\,\widetilde{q}(E)\,\,(\not= \widetilde{V})$ is 
$\widetilde{f}(\widetilde{\frak g})$-invariant  
by means of $(3.4)$.  
This contradicts the irreducibility of 
$(\widetilde{f},\widetilde{q},\widetilde{V})$. 
So we have $\widetilde{q}(E)=0$. 
However, that 
$\widetilde{q}(E)=0$ contradicts the assumption 
that Ker$(\widetilde{q})={\frak k}$. 
So, we have $a\not=0$.  
\par
Notice here that, if 
we put $E'=\frac{1}{a}\,E$, 
then, we have also
$\widetilde{\frak g}={\frak g}\oplus{\mathbb R}E'$, and 
$[E',\widetilde{\frak g}]=\{0\}$. 
Furthermore, the affine representation $(\widetilde{f},\widetilde{q},\widetilde{V})$ of $\widetilde{\frak g}$ 
still satisfies all the conditions $(1)$, $(2)$ and $(3)$ 
for the set 
${\mathcal PF}(\widetilde{\frak g},{\frak k})$ in Definition 3.1. Because, for $(1)$, $(2)$ they are the same, 
and for $(3)$, we have that 
$\widetilde{f}(E')=\frac{1}{a}\,\widetilde{f}(E)=I$, and 
$\widetilde{q}(E')=\frac{1}{a}\,\widetilde{q}(E)\not=0$. 
Thus, $(\widetilde{f},\widetilde{q},\widetilde{V})\in {\mathcal PF}_0(\widetilde{\frak g},{\frak k})$. 
\vskip0.3cm
\par
\underline{Case $(b)$}: 
$V_1\oplus\cdots\oplus V_s=\{0\}$. 
\par
In this case, 
let $\lambda$ be any non-zero 
$\lambda^i_j$ in $(ii)$, and consider a non-zero complex subspace 
$$W:=\sum_{\lambda^k_{\ell}=\lambda}
({\mathbb C}v^k_{\ell}\oplus
{\mathbb C}\overline{v^k_{\ell}})
$$
of $\sum_{j=s+1}^r{V_j}^{\mathbb C}$, 
where 
$\lambda^k_{\ell}$ run over the set of all 
complex eigenvalues of 
$\widetilde{f}(E)^{\mathbb C}$ in (3.1)  which 
are equal to $\lambda$. 
Then, 
\begin{equation}
W\cap \widetilde{V}=
\sum_{\lambda^k_{\ell}=\lambda}
({\mathbb R}u^k_{\ell}\oplus {\mathbb R}w^k_{\ell}),
\end{equation}
where we denote $v^k_{\ell}=u^k_{\ell}+\sqrt{-1}w^k_{\ell}$, 
$u^k_{\ell}, \,w^k_{\ell}\in \widetilde{V}$. 
\par 
Then, for each $\widetilde{X}\in \widetilde{\frak g}$, 
$$
\widetilde{f}(E)^{\mathbb C}\widetilde{f}(\widetilde{X})^{\mathbb C}v^k_{\ell}=
\widetilde{f}(\widetilde{X})^{\mathbb C}
\widetilde{f}(E)^{\mathbb C}v^k_{\ell}
=\lambda\,\widetilde{f}(\widetilde{X})^{\mathbb C}v^k_{\ell},
$$
which implies that $\widetilde{f}(\widetilde{X})^{\mathbb C} v^k_{\ell}$ belongs to $W$. 
Furthermore, we have for each $\widetilde{X}\in \widetilde{\frak g}$, 
$$
\widetilde{f}(\widetilde{X})(W\cap \widetilde{V})\subset
W\cap \widetilde{V}. 
$$
Because, since  it holds that 
$$
\widetilde{f}(\widetilde{X})u^k_{\ell}
+\sqrt{-1}\widetilde{f}(\widetilde{X})w^k_{\ell}
=\widetilde{f}(\widetilde{X})^{\mathbb C}v^k_{\ell}
=\sum_{\lambda^p_q=\lambda}
(\alpha_{pq}v^p_q+\beta_{pq}\overline{v^p_q}),
$$
for some complex numbers 
$\alpha_{pq}$, and $\beta_{pq}$, 
we have 
\begin{align*}
\widetilde{f}(\widetilde{X})u^k_{\ell}&=
\sum_{\lambda^p_q=\lambda}
\left\{
\Re e(\alpha_{pq}+\beta_{pq})\,u^p_q
+\Im m(-\alpha_{pq}+\beta_{pq})\,w^p_q
\right\}\in W\cap \widetilde{V},\\
\widetilde{f}(\widetilde{X})w^k_{\ell}&=
\sum_{\lambda^p_q=\lambda}
\left\{
\Im m(\alpha_{pq}+\beta_{pq})\,u^p_q
+\Re e(\alpha_{pq}-\beta_{pq})\,w^p_q
\right\}\in W\cap \widetilde{V}.
\end{align*}
Thus, together with (3.5), we have 
$\widetilde{f}(\widetilde{X})(W\cap \widetilde{V})\subset W\cap \widetilde{V}$. 
\par
Since $(\widetilde{f},\widetilde{q},\widetilde{V})$ is an irreducible representation of $\widetilde{\frak g}$, 
we have 
$$
\widetilde{V}=W\cap \widetilde{V}=
\sum_{\lambda^k_{\ell}=\lambda}
({\mathbb R}u^k_{\ell}+{\mathbb R}w^k_{\ell}),
$$
which means that 
$V_{s+1}{}^{\mathbb C}\oplus\cdots\oplus V_r{}^{\mathbb C}$ is the sum of the two eigenspaces of $\widetilde{f}(E)^{\mathbb C}$ 
with the eigenvalues $\lambda$ and $\overline{\lambda}$  
for some complex number 
$\lambda=a+\sqrt{-1}b$ $(a,b\in {\mathbb R})$ with 
$b\not=0$. 
\par
By means of (3.2), 
$\widetilde{f}(E)$ can be written as 
\begin{equation}
\widetilde{f}(E)=a\,I+b\,J, 
\end{equation}
where $I$ is the identity transformation of $\widetilde{V}$, and 
$J$ is the transformation of $\widetilde{V}$ of the form: 
\begin{equation}
J(u_j)=-w_j,\quad J(w_j)=u_j \quad (j=1,\dots,d), 
\end{equation}
where $\{u_1,\dots,u_d,w_1,\dots,w_d\}$ is a basis of 
$\widetilde{V}$ $(\dim \widetilde{V}=2d)$. 
Then, $J$ is a complex structure of $\widetilde{V}$, i.e., 
$J^2=-I$. 
\par
The complex structure $J$ of $\widetilde{V}$ is 
$\widetilde{f}(\widetilde{\frak g})$-invariant, 
i.e., 
\begin{equation}
J\,\widetilde{f}(\widetilde{X})=\widetilde{f}(\widetilde{X})\,J
\quad (\forall\,\,\widetilde{X}\in \widetilde{\frak g}). 
\end{equation}
Because, since 
$E$ is central in $\widetilde{\frak g}$, we have, 
for each $\widetilde{X}\in \widetilde{\frak g}$, 
\begin{equation}
\widetilde{f}(E)\widetilde{f}(\widetilde{X})
=\widetilde{f}(\widetilde{X})\widetilde{f}(E).
\end{equation}
The left hand side of $(3.9)$ coincides with 
$$
(a\,I+b\,J)\,\widetilde{f}(\widetilde{X})= 
a\,\widetilde{f}(\widetilde{X})+b\,J\,\widetilde{f}(\widetilde{X}). 
$$
The right hand side of $(3.9)$ is equal to 
$$
\widetilde{f}(\widetilde{X})(a\,I+b\,J)=a\,\widetilde{f}(\widetilde{X})
+b\,\widetilde{f}(\widetilde{X})\,J. 
$$
Since $b\not=0$, we have (3.8). 
\par 
Notice that 
$\widetilde{q}(E)\not=0$ since Ker$(\widetilde{q})={\frak k}$. 
Since 
$$
\widetilde{q}([\widetilde{X},\widetilde{Y}])=
\widetilde{f}(\widetilde{X})(\widetilde{q}(\widetilde{Y}))
-\widetilde{f}(\widetilde{Y})(\widetilde{q}(\widetilde{X})) 
\quad (\widetilde{X},\widetilde{Y}\in \widetilde{\frak g}),
$$
we have 
\begin{equation}
\widetilde{f}(\widetilde{Y})(\widetilde{q}(E))
=\widetilde{f}(E)(\widetilde{q}(\widetilde{Y}))
=(a\,I+b\,J)(\widetilde{q}(\widetilde{Y})).
\end{equation}
Then, we have, for all $\widetilde{Y}\in \widetilde{\frak g}$, 
\begin{align}
\widetilde{q}(\widetilde{Y})
&=\frac{1}{a^2+b^2}(a\,I-b\,J)
\left(\widetilde{f}(\widetilde{Y})(\widetilde{q}(E))\right)\nonumber\\
&=\widetilde{f}(\widetilde{Y})\left(\frac{1}{a^2+b^2}(a\,I-b\,J)
(\widetilde{q}(E))\right)\nonumber\\
&=\widetilde{f}(\widetilde{Y})v_0
\end{align}
where 
\begin{equation}
0\not=v_0:=
\frac{1}{a^2+b^2}(a\,I-b\,J)(\widetilde{q}(E))\in \widetilde{V}.
\end{equation}
By (3.11), $\widetilde{f}({\frak k})v_0=\{0\}$ 
since Ker$(\widetilde{q})={\frak k}$. 
Since 
$\widetilde{q}:\,\widetilde{\frak g}\rightarrow \widetilde{V}$ is a surjection 
and $\widetilde{\frak g}={\frak g}\oplus{\Bbb R}E$, and (3.11), 
we also have 
\begin{equation}
\widetilde{V}=\widetilde{f}({\frak g})v_0\oplus{\Bbb R}\,\widetilde{f}(E)v_0.
\end{equation}
Thus, we have 
$(\widetilde{f},\widetilde{q},\widetilde{V})
\in {\mathcal F}^{I\!I}_0(\widetilde{\frak g},{\frak k})$. 
\end{pf}
\vskip0.6cm\par
\begin{rem} \,\,$(1)$ \,\, 
The set ${\mathcal PF}_0({\widetilde{\frak g}},{\frak k})$ 
corresponds 
$($\cite{U}, Theorem 1.3$)$ to the set of all real irreducible representations 
$(f,\widetilde{V})$ of ${\frak g}$ of dimension $\dim M+1$ which admits 
a nonzero vector $v_0\in \widetilde{V}$ such that 
$f({\frak k})v_0=\{0\}$ and $\widetilde{V}=f({\frak g})v_0\oplus 
{\Bbb R}v_0$. 
Each $(f,\widetilde{V})$ induces a $G$-equivariant
centro-affine immersion 
$\varphi:\, M=G/K\rightarrow \widetilde{V}$ by 
$\varphi(xK)=f(xK)v_0$ $(xK\in G/K)$ with transversal vector field 
$\xi_{xK} =-\overrightarrow{o\,\varphi(xK)}$ $(xK\in G/K)$ 
$($\cite{NS}$)$.  
Our Theorem 3.3 suggests that the set 
${\mathcal F}^{I\!I}_0(\widetilde{\frak g},{\frak k})$ 
would correspond to the set 
of all $G$-equivariant affine immersions given by 
$\varphi'(xK)=\widetilde{f}(xK)v_0$ $(xK\in G/K)$ with the transversal vector field 
$\xi'_{xK} =
\overrightarrow{\varphi'(xK)\,\widetilde{f}(E)\varphi'(xK)}$ $(xK\in G/K)$. 
Then, $M=G/K$ would admit a $G$-invariant affine connection via 
the immersion $\varphi'$ $($cf. \cite{NS}$)$. 
\par
$(2)$ \,\, There is an example belonging to  
${\mathcal F}^{I\!I}_0(\widetilde{\frak g}, {\frak k})$. 
Indeed, 
let us recall Example 11.1 in the book of Takeuchi \cite{T}, p. 119, and 
let
${\frak g}={\frak su}(2)$, ${\frak k}=\{0\}$, 
$\widetilde{\frak g}={\frak g}\oplus {\Bbb R}E$, with $E=I_2$. 
Let 
$V={\Bbb C}^2
=\left\{\begin{pmatrix}
z\\w
\end{pmatrix};\,z,w\in {\Bbb C}\right\}$, 
and 
$\widetilde{V}$, the real 4-dimensional space $V$ 
restricted to the field ${\Bbb R}$, and 
${\frak g}$ acts on $\widetilde{V}$ by the matrix multiplications
$\widetilde{f}(X)v=Xv\in \widetilde{V}$, $(X\in {\frak g},v\in \widetilde{V})$,
and 
$\widetilde{f}(E)v=I_2v=v$ $(v\in \widetilde{V})$. Then, 
$\widetilde{V}$ admits the complex structure 
$J$ defined by 
$$
Jv=\begin{pmatrix}
0&-1\\
1&0
\end{pmatrix}\overline{v},\quad v\in \widetilde{V}, 
$$
which satisfies $J(Xv)=XJv$, $(X\in {\frak g},v\in \widetilde{V})$. 
The linear mapping 
$\widetilde{q}:\widetilde{\frak g}\rightarrow \widetilde{V}$ is given by 
$\widetilde{q}\left(
\begin{pmatrix}
i\theta&\alpha\\
-\overline{\alpha}&-i\theta
\end{pmatrix}+\xi E
\right)
=\begin{pmatrix}
i\theta\\
-\overline{\alpha}
\end{pmatrix}
+\xi\begin{pmatrix}
1\\0
\end{pmatrix}, 
$
and satisfies 
$$
\widetilde{q}([\widetilde{X},\widetilde{Y}])
=\widetilde{f}(\widetilde{X})\widetilde{q}(\widetilde{Y})
-\widetilde{f}(\widetilde{Y})\widetilde{q}(\widetilde{X}),\quad 
(\widetilde{X},\widetilde{Y}\in \widetilde{\frak g}).
$$
Then, 
$(\widetilde{f},\widetilde{q},\widetilde{V})\in 
{\mathcal F}^{I\!I}_0(\widetilde{\frak g},{\frak k})$. 
\par
$(3)$\,\, The real representations in 
${\mathcal F}^{I\!I}_0(\widetilde{\frak g},{\frak k})$ 
in the case ${\frak k}=\{0\}$ were treated, called 
the representations of $\widetilde{\frak g}$ {\em of class II}, 
in the book of Takeuchi $($\cite{T}, pp. 85--92$)$. 
\par
$(4)$ The union in the right hand side for 
${\mathcal F}_0(\widetilde{\frak g},{\frak k})$ 
in Theorem 3.3 seems to be a disjoint union. 
\end{rem}
\vskip2cm\par       


\begin{thebibliography}{99}
\bibitem{A} Y. Agaoka, {\em Invariant flat projective structures on homogeneous spaces}, \newblock Hokkaido Math. J., {\bf 11} (1982), \newblock 125--172. 
\bibitem{D} H. Doi, {\em Non-existence of torsion free flat connections on reductive homogeneous spaces}, Hiroshima Math. J., {\bf 9} (1979), 321--322.
\bibitem{E} A. Elduque, \newblock 
{\em Invariant projectively flat affine connections on Lie groups}, \newblock
Hokkaido Math. J., {\bf 30} (2001), \newblock 231--239. 
\bibitem{H} S. Helgason, \newblock 
{\em Differential Geometry, Lie Groups, and Symmetric Spaces}, 
Academic Press, New York, 1978. 
\bibitem{KN} S. Kobayashi and K. Nomizu, \newblock 
{\em Foundation of Differential Geometry, Vol. I, II}, (1963), (1969), 
John Wiley and Sons, New York. 
\bibitem{MO} H. Matsushima and K. Okamoto, \newblock 
{\em Non-existence of torsion free flat connections on a real semisimple Lie group}, \newblock Hiroshima Math. J., {\bf 9} (1979), \newblock 59--60. 
\bibitem{MM} Y. Matsushima and S. Murakami, \newblock 
{\em On vector bundle valued harmonic forms and automorphic forms on symmetric Riemannian manifolds}, \newblock Ann. Math., {\bf 78} \newblock (1963), 365--416. 
\bibitem{N} K. Nomizu, \newblock {\em Invariant affine connections on homogeneous spaces}, 
Amer. J. Math., {\bf 76}\newblock (1954), 33--65. 
\bibitem{NS} K. Nomizu and T. Sasaki, \newblock 
{\em Affine Differential Geometry, Geometry of Affine Immersions}, 
(1994), Shokabo, Tokyo; (1994), Cambridge Univ. Press, Cambridge.
\bibitem{S1} H. Shima, \newblock 
{\em The Geometry of Hessian  Structures}, (2001), Shokabo, Tokyo; (2007), World Sci. Publ., Hackensack.  
\bibitem{S2} \bysame, \newblock 
{\em Homogeneous spaces with invariant projectively flat affine connections}, \newblock 
Trans. Amer. Math. Soc., {\bf 351}\newblock 
(1999), \newblock 4713--4726. 
\bibitem{T} M. Takeuchi, \newblock 
{\em Real Irreducible Representations of Real Simple Lie Algebras (in Japanese)}, 
Lecture Notes Series, Vol. 3,  
\newblock 1996, \newblock Osaka University. 
\bibitem{U} H. Urakawa, \newblock 
{\em On invariant projectively flat affine connections}, \newblock 
Hokkaido Math. J., {\bf 28} \newblock (1999), \newblock 333--356. 
\end{thebibliography}
\end{document}